\documentclass[12pt]{article}

\usepackage[T2A]{fontenc}
\usepackage[utf8]{inputenc}

\usepackage{amsfonts}
\usepackage{amssymb}
\usepackage{amsmath}
\usepackage{amsthm}

\def \le {\leqslant}
\def \ge {\geqslant}

\begin{document}

\begin{Huge}
 \centerline{\bf On Diophantine exponents in dimension 4}
 
\end{Huge}
\vskip+0.5cm

\centerline{ by  {\bf Dmitry Gayfulin} and {\bf Nikolay Moshchevitin}\footnote{ Research is supported by RFBR grant No.12-01-00681-a and by the grant of Russian Government, project 11.  
G34.31.0053 and by the grant NSh-2519.2012.}}
\vskip+0.5cm

{\bf 1. Introduction.}

We consider a vector $\Theta = (\theta_1,...,\theta_n), n\ge 2$ and suppose that the numbers $1,\theta_1,...,\theta_n$ are linearly independent over $\mathbb{Z}$.
Put
$$
\psi_\Theta (t) =
\min_{ q \in \mathbb{Z}_+,  q\le t}
\, \max_{1\le j \le n} ||q\theta_j||.
 $$
We consider the  ordinary Diophantine exponent
$\omega = \omega(\Theta)$ and the uniform Diophantine exponent $\hat{\omega} = \hat{\omega} (\Theta)$ defined as 
$$
\omega = \omega(\Theta)=
\sup \left\{
\gamma:\,\,\liminf_{t\to+\infty}
t^\gamma \psi_\Theta (t) 
<+\infty\right\},
$$
$$
\hat{\omega} = \hat{\omega}(\Theta)=
\sup \left\{
\gamma:\,\,\limsup_{t\to+\infty}
t^\gamma \psi_\Theta (t) 
<+\infty\right\}.
$$ 
It is clear that
$$
\frac{1}{n} \le \hat{\omega}\le 1
$$
and
\begin{equation}\label{tree}
 \omega \ge \hat{\omega}.
\end{equation}
In \cite{J} 
 V. Jarn\'{\i}k proved that in the case $ n = 2$ the trivial inequality
(\ref{tree}) may be improved to 
\begin{equation}\label{simult}
 \omega \ge  \frac{\hat{\omega}^2}{1-\hat{\omega}}.
\end{equation} 
In \cite{Lo}
M. Laurent proved that the bound (\ref{simult})
is optimal for the case $n=2$.
In \cite{MCZ} the author in the case $n=3$ proved the inequality

\begin{equation}\label{shmu}
\omega  \ge \frac{\hat{\omega}}{2} \left( \frac{\hat{\omega}}{1-\hat{\omega}}
+\sqrt{\left(\frac{\hat{\omega}}{1-\hat{\omega}}\right)^2
+\frac{4\hat{\omega}}{1-\hat{\omega}}}\right) .
\end{equation}
A different proof of the inequality (\ref{shmu})
was given by W. Schmidt and L. Summerer \cite{SSu} by means of a new powerful method developed in \cite{ssu,ssu2}.

It turned out that the inequality (\ref{shmu}) is optimal. In \cite{SSu} it was announced that D. Roy recently obtained such a result (see also author's announcement from Section 3.3 from \cite{arno} concerning Theorem 14.)
However  no optimal inequality is 
known in dimensions $ n \ge 4$.
Probably the best known general inequality is due to W.M. Schmidt and L. Summerer.  It is as follows. For an arbitrary $n\ge 2$ one has
$$
\omega  \ge
\frac{\hat{\omega}^2 + (n-2)\hat{\omega} }{(n-1)(1-\hat{\omega})}
.
$$
In the present paper we obtain some results concerning the case $ n = 4$.

{\bf 2. Index $\hbox{\got i}(\Theta)$.}
\,\,
We consider a vector $\Theta = (\theta_1,\theta_2,\theta_3,\theta_4)$ such that real numbers $1,\theta_1,\theta_2,\theta_3,\theta_4$
are linearly independent over $\mathbb{Z}$.
Here we introduce the value  $\hbox{\got i}(\Theta)$ which is of importance for the formulation of our results. We define it to be the {\it index} of the vector $\Theta$.
We consider the sequence ${\bf z}_\nu =(q_\nu, a_{1,\nu},a_{2,\nu} a_{3,\nu},a_{4,\nu})\in \mathbb{Z}^5,\,\, \nu =0,1,2,3,4,...$
of all best approximations to the vector $\Theta$, so that
$$
q_0 <q_1<....<q_{\nu}<q_{\nu+1}<...,\,\,\,\,\,\,  \zeta_0>\zeta_1>...>\zeta_\nu>\zeta_{\nu+1}>... ,
$$
where
$$
\zeta_\nu = \max_{1\le j \le 4} ||q_\nu\theta_j||,\,\,\,\, ||q_\nu \theta_j || = |q_\nu \theta_j - a_{\nu,j}|
$$
and
$$
\zeta_\nu =
\min_{x\in \mathbb{Z},\, 1\le x \le q_{\nu}}\,\,  \max_{1\le j \le 4} ||x\theta_j||
.
$$
Take $\alpha < \hat{\omega}(\Theta )$. Then
\begin{equation}\label{aal}
\zeta_\nu\le q_{\nu+1}^{-\alpha}
\end{equation}
for $\nu$ large enough.

It is a well known fact (see \cite{M}, Section 4.1) that for any $\nu_0$ all the vectors ${\bf z}_\nu, \nu \ge \nu_0$ cannot lie in a common linear subspace $L \subset \mathbb{R}^5$
of dimension ${\rm dim }\, L \le 4$. V. Jarn\'{\i}k \cite{J} proved that there exists infinitely many $\nu$ such that three vectors ${\bf z}_{\nu-1}, {\bf z}_\nu, {\bf z}_{\nu+1}$
are linearly independent. One can easily deduce from these two  facts the following

{\bf Proposition 1.}\,\, {\it There exist infinitely many sets of indices $ (\nu=r_0,r_1,r_2,\ldots,r_n,r_{n+1}=k)_l$ 
(here  $ n =n_l$  depends on $l$) such that\footnote{We have in mind that all the parameters are depend on $l$. However, sometimes we will not use the lower index ${}_l$ with a view to simplify the notation}

{\rm (i)}  \,\, $\nu \to \infty $ as $ l\to \infty$;

{\rm (ii)} \,\, The triplets of vectors  $({\bf z}_{r_i-1}, {\bf z}_{r_i}, {\bf z}_{r_i+1})$ are linearly independent for any $0\le i\le n+1$;

{\rm (iii)} for any $ i$ from the interval  $1\le i \le n$ each two consecutive couples of vectors 
$$
({\bf z}_{r_{i}}, {\bf z}_{r_i+1}),
({\bf z}_{r_{i}+1}, {\bf z}_{r_i+2})
,...,
({\bf z}_{r_{i+1}-1}, {\bf z}_{r_{i+1}})
$$ 
lie in the same two-dimencsional linear subspece; we denote this subspace by $L_{i}$, so
$$
L_{i} =
{\rm span } \, ( {\bf z}_{r_{i}}, {\bf z}_{r_{l+1}}) =
{\rm span } \,( {\bf z}_{r_{i+1}-1}, {\bf z}_{r_{i+1}}  );
$$

{\rm (iv)}\,\, there exists a three-dimensional linear  subspace  $T_l$ such that 
for any $ i$ from the interval $1\le i \le n_l$ the triple $({\bf z}_{r_i-1}, {\bf z}_{r_i}, {\bf z}_{r_i+1})$
belongs to this subspace, so
$$
T_l = {\rm span}\,({\bf z}_{r_i-1}, {\bf z}_{r_i}, {\bf z}_{r_i+1}),\,\,\,\,
1\le i \le n_l;
$$

 {\rm (v)} \,\,  ${\bf z}_{\nu_l -1}\not\in T_l$  and ${\bf z}_{k_l +1}\not\in T_l$;

 {\rm (vi)} \,\, the collection $ {\bf z}_{\nu-1}, {\bf z}_{ r_{n}-1}, {\bf z}_{k-1},  {\bf z}_{k},  {\bf z}_{k+1}$
consists of five linearly independent vectors.}

It may happen that for the vector $\Theta$ there exists an integer $n$  and a sequence of indices $ (\nu,r_1,r_2,\ldots,r_n,k)_l$ such that such that the conclusions
{\rm (i) - (vi)} of Proposition 1 hold for $n_l =n$. Then we define
\begin{equation}
\begin{split}
\hbox{\got i}(\Theta) = \min \{n :\, \text{there exists a sequence } (\nu,r_1,r_2,\ldots,r_n,k)_l\\ \text{ such that the conditions
{\rm (i) - (vi)}   hold for}\,  n_l =n \}.
\end{split}
\end{equation}
If it is not so we define $\hbox{\got i}(\Theta) = \infty$.

For example  in the case  when there exist infinitely many $\nu$ such that  every five vectors
$${\bf z}_{\nu-1}, {\bf z}_{\nu}, {\bf z}_{\nu+1}, {\bf z}_{\nu+2}, {\bf z}_{\nu+3}$$
 are linearly independent,  one has $ \hbox{\got i}(\Theta)=1$. In fact, 
certain lower bounds for $\omega(\Theta )$ in terms of $\hat{\omega}(\Theta )$ in this case were discussed in  \cite{arno}, Subsection 3.5, footnote $^7$. 
We believe that  for simultaneous approximation these bounds are optimal in the case  $ \hbox{\got i}(\Theta)=1, 2$.

One can easily see that  $\hbox{\got i}(\Theta)$ can attain any positive integer value for a vector $\Theta$ with the components
$\theta_1, \theta_2,\theta_3,\theta_1$,
linearly independent over $\mathbb{Z}$ together with 1.

{\bf 3. The main result.}

We consider the  polynomial
$$
f_1(x)  =f_{1}^{[\hat{\omega}]} (x)=
x^3 -\frac{\hat{\omega}}{1-\hat{\omega}} x^2 -
\frac{\hat{\omega}}{1-\hat{\omega}} x-
\frac{\hat{\omega}}{1-\hat{\omega}} .
$$
For every $\hat{\omega} \in \left[\frac{1}{4},1\right)$ it has a unique real positive root
 $G_1 (\hat{\omega})$.
One can see that   $G_1 \left(\frac{1}{4}\right)=1$ and  $G_1 (\hat{\omega})$ increases to infinity as $\hat{\omega}$ increases to 1.

Then we consider two polynomials
$$
f_{2,1} (x)  =f_{2,1}^{[\hat{\omega}]} (x) =x^4 -\frac{\hat{\omega}}{1-\hat{\omega}} x^3 -
\frac{\hat{\omega}}{1-\hat{\omega}} x^2+
\left(\frac{\hat{\omega}}{1-\hat{\omega}}\right)^2 x-
\frac{\hat{\omega}}{(1-\hat{\omega})^2} 
$$
and
$$
f_{2,2}(x)   =f_{2,2}^{[\hat{\omega}]} (x) =
x^4 -\frac{\hat{\omega}}{1-\hat{\omega}} x^3 -
\frac{\hat{\omega}}{1-\hat{\omega}} x^2+
\frac{\hat{\omega}}{1-\hat{\omega}} x-
\frac{\hat{\omega}}{(1-\hat{\omega})^2}
.$$
One can see that for $\hat{\omega} = \frac{1}{2}$ one has $f_{2,1} = f_{2,2}$.
For every $\hat{\omega} \in \left[\frac{1}{4},\frac{1}{2}\right]$ the polynomial $f_{2,1}$ has the unique positive root
  $G_{2,1} (\hat{\omega})$. 
For every $\hat{\omega} \in \left[\frac{1}{2},1\right)$
the polynomial $f_{2,2} $ has the unique real positive root 
  $G_{2,2} (\hat{\omega})$.
  Put
$$
G_{2} (\hat{\omega}) = \max (G_{2,1} (\hat{\omega}),G_{2,2} (\hat{\omega}))
=
\begin{cases}
G_{2,1} (\hat{\omega}) \,\,\,\text{if}\,\,\, \frac{1}{4} \le \hat{\omega} \le \frac{1}{2},\cr
G_{2,2} (\hat{\omega}) \,\,\,\text{if}\,\,\, \frac{1}{2} \le \hat{\omega}<1,
\end{cases}
$$
It can be easily seen that    $G_2 \left(\frac{1}{4}\right)=1$ and  $G_2 (\hat{\omega})$ increases to infinity as $\hat{\omega}$ increases to 1.
It is clear that 
\begin{equation}\label{expla}
G_2(\hat{\omega}) <G_1(\hat{\omega}) \,\,\,\,\,\text{for}\,\,\,\,\, \frac{1}{4}<\hat{\omega}<1.
\end{equation}
We also consider the polynomial
\begin{equation}\label{pool}
f_3(x)  =f_{3}^{[\hat{\omega}]} (x)=
x^5 -\frac{\hat{\omega}}{1-\hat{\omega}} x^4 -
\frac{\hat{\omega}}{1-\hat{\omega}} x^3+\frac{{\hat{\omega}}^2}{({1-\hat{\omega}})^2}x^2-
\frac{\hat{\omega}}{({1-\hat{\omega}})^2} .
\end{equation}
Denote $G_3(\hat{\omega})$ a unique positive root of $f_3(x)$. We can easily see that $G_3(\frac14)=1,G_3(\hat{\omega}) <G_2(\hat{\omega})$ for $\frac{1}{4}<\hat{\omega}<1$  and  $G_3 (\hat{\omega})$ increases to infinity as $\hat{\omega}$ increases to 1.

{\bf Theorem 1.}\,\,{\it

{\rm 1)} Suppose that  $\hbox{\got i}(\Theta) = 1$. Then

\begin{equation}\label{oo1}
\omega (\Theta) \ge \hat{\omega} (\Theta) G_1(\hat{\omega}(\Theta)).
\end{equation}

{\rm 2)} Suppose that $\hbox{\got i}(\Theta) = 2.$   Then
\begin{equation}\label{oo2}
\omega (\Theta) \ge \hat{\omega} (\Theta) G_2(\hat{\omega}(\Theta)).
\end{equation}

{\rm 3)} Suppose that $\hbox{\got i}(\Theta) = 3$.   Then
\begin{equation}\label{oo3}
\omega (\Theta) \ge \hat{\omega} (\Theta) G_3(\hat{\omega}(\Theta)).
\end{equation}
}

 It is not difficult to show thet in the case $\hbox{\got i}(\Theta) = 1$ Theorem 1 gives the opimal lower bound. However we do not prove this fact in the present paper.
The bound (\ref{oo2}) of the case 
$\hbox{\got i}(\Theta) = 2$  
should also be optimal. Moreover, we believe that there exists a constant $\hat{\omega_0}$ such that the bound (\ref{oo3}) is optimal in the interval $\frac14\le\hat{\omega}\le\hat{\omega_0}$.

In the sequel for the case $ \hbox{\got i}(\Theta) = 1$ we give a complete proof of the statement, as well as for the case 
 $ \hbox{\got i}(\Theta) = 2$ and $\hat{\omega}\le \frac{1}{2}$.
For the cases 
 $ \hbox{\got i}(\Theta) = 2$ and $\hat{\omega}\ge  \frac{1}{2}$ and
 $ \hbox{\got i}(\Theta) = 3$.

{\bf Remark 1.} In any cases Proposition 1 gives a system of equalities and inequalities.
Some of these inequalities are essential for the case under the consideration, but some of them are not.
We performed computer calculations  to determine which inequalities are essential for each case. In such a way we 
found the systems (\ref{zis},\ref{zis2},\ref{zis3})
below. As these systems
define simplicial cones in the considered spaces, this gives  proofs of the main estimates.

{\bf Remark 2.}
We performed more extensive calculations in the case  $ \hbox{\got i}(\Theta) = 3$ and found out that sometimes  for
$\hat{\omega}>\hat{\omega_0}$
 the polynimial $f_3$ from (\ref{pool}) does not give an optimal lower bound in (\ref{oo3}).
For some values of $\hat{\omega}$ better bounds may be obtained due to the polynomials
$$
x^5-x^4\frac{\hat{\omega}}{1-\hat{\omega}}-x^3\frac{\hat{\omega}}{1-\hat{\omega}}-\hat{\omega}\frac{(1-x\hat{\omega})^2}{(1-\hat{\omega})^3},
\,\,\,\,\,\,
x^5-x^4\frac{\hat{\omega}}{1-\hat{\omega}}-x^3\frac{\hat{\omega}}{1-\hat{\omega}}-\hat{\omega}\frac{(1-x\hat{\omega})(2-x)}{(1-\hat{\omega})^2}.
$$
However the calculations are too cumbersome and do not rely on any new ideas.

{\bf 4. Case $ \hbox{\got i}(\Theta) = 1$.}\,\, Here one must note that if
there are two indices $j_1\le j_2$ such that
$${\rm span } \, (  {\bf z}_{j_1}, {\bf z}_{j_1+1}) =
{\rm span } \,(  {\bf z}_{j_2-1}, {\bf z}_{j_2})
$$
then
\begin{equation}\label{twailait}
\zeta_{j_1}q_{j_1+1}\asymp_\Theta \zeta_{j_2-1} q_{j_2}.
\end{equation}
This is a well known statement (see Lemma 1 from \cite{MCZ} or Theorem 2.13  from \cite{CHEU}).
Now from (vi) we see that $5\times 5$  determinant constructed from the coordinates of vectors
$$
{\bf z}_{\nu-1}, {\bf z}_{r_1-1},{\bf z}_{k-1},{\bf z}_k,{\bf z}_{k+1}
$$
is $ \ge 1$ in absolute value. So 
\begin{equation}\label{oole}
1\ll_\Theta \zeta_{\nu-1}\zeta_{r_1-1}\zeta_{k-1}\zeta_{k}q_{k+1}.
\end{equation} 
We have
\begin{equation}\label{oole1}
\zeta_\nu q_{\nu+1}\asymp_\Theta \zeta_{r_1-1} q_{r_1},\,\,\,\,
\zeta_{r_1}q_{r_1+1}\asymp_\theta \zeta_{k-1}q_k.W
\end{equation}
Now we take two parameters $u,v \in (0,1) $ defined  from the system of equalities
$$
\frac{\alpha}{(1-\alpha)u} =-v +\frac{1}{1-\alpha} =\frac{(1-\alpha)u+\alpha}{(1-\alpha)(1-v)}.
$$
From (\ref{oole}) we have
$$
\zeta_{\nu-1}\zeta_\nu q_{\nu+1} \cdot
\zeta_{r_1-1}\zeta_{r_1}q_{r_1+1}
\cdot
\zeta_{k-1}\zeta_{k}q_{k+1}
\gg_\Theta
\zeta_\nu q_{\nu+1} \cdot\zeta_{r_1}q_{r_1+1}
\asymp_\Theta
(\zeta_\nu q_{\nu+1})^{1-u}
\cdot
(\zeta_{r_1-1}q_{r_1})^{u}
(\zeta_{r_1}q_{r_1+1})^{v}
\cdot
( \zeta_{k-1}q_k.)^{1-v}.
$$
So at least one of the following three inequalities is valid:
\begin{equation}\label{odin}
\zeta_{\nu-1}\zeta_\nu q_{\nu+1}\gg_\Theta (\zeta_\nu q_{\nu+1})^{1-u}
,\end{equation}
\begin{equation}\label{dva}
\zeta_{r_1-1}\zeta_{r_1}q_{r_1+1}\gg_\Theta (\zeta_{r_1-1}q_{r_1})^{u}
(\zeta_{r_1}q_{r_1+1})^{v}
,\end{equation}
\begin{equation}\label{tri}
\zeta_{k-1}\zeta_{k}q_{k+1}
\gg_\Theta ( \zeta_{k-1}q_k.)^{1-v}
.\end{equation}
We  take $\alpha < \hat{\omega}(\Theta )$ close to $\hat{\omega}$. Then by  (\ref{aal}) and the choice of parameters $u,v$ we deduce te following.
Given any small positive  $\varepsilon$,
if (\ref{odin}) holds then $ q_{\nu+1}\gg_\Theta q_{\nu}^{G_1(\hat{\omega})-\varepsilon}$,
if (\ref{dva}) holds then $ q_{r_1+1}\gg_\Theta q_{r_1}^{G_1(\hat{\omega})-\varepsilon}$,
if (\ref{tri}) holds then $ q_{k+1}\gg_\Theta q_{k}^{G_1(\hat{\omega})-\varepsilon}$.

So in the case $ \hbox{\got i}(\Theta) = 1$
everything is proven.

{\bf 5. Case  $ \hbox{\got i}(\Theta) = 2$ and $\hat{\omega}\le \frac{1}{2}$.}

Put
$$
\hbox{\got X} = (\xi_{\nu-1},\xi_\nu, \xi_{r_1-1},\xi_{r_1},\xi_{r_2-1},\xi_{r_2},\xi_{k-1},\xi_k, X_\nu, X_{\nu+1},
X_{r_1}, X_{r_1+1}, X_{r_2}, X_{r_2+1}, X_k, X_{k+1}).
$$
Then $\hbox{\got X}$ is a point in $16$-dimensional space $ {\cal R} =\mathbb{R}^{16}$.
We consider  subspace 
$${\cal R}' = \{ \hbox{\got X} \in {\cal R}:\,\,\, X_\nu = 0\}.$$
We are interested in the system  in fifteen inequalities
\begin{equation}\label{zis}
\begin{cases}
\xi_{\nu-1}+\xi_{r_2-1}+\xi_{k-1}+\xi_k+ X_{k+1} \ge 0,\cr
\xi_{j-1}-\alpha X_j \le 0,\,\,\, j = \nu,\nu+1,r_1+1,r_2+1, k+1,\cr
X_{r_1+1}\le X_{r_2},\,\,\, X_{r_2+1} \le X_k,\cr
\xi_{r_1-1}\ge \xi_{r_1}\ge\xi_{r_2-1},\,\,\, \xi_{r_2}\ge \xi_{k-1},\cr
X_{j+1} \le gX_{j},\,\,\, j =\nu, r_1,r_2, k.  
\end{cases}
\end{equation}
and in the equation
\begin{equation}\label{eq1}
\xi_\nu +X_{\nu+1} =\xi_{r_1-1}+X_{r_1}
\end{equation}

{\bf Lemma 1.}\,\,{\it Suppose that $\alpha, g >0$.
Suppose that there exists $\hbox{\got X} \in {\cal R}'$ such that its coordinates satisfy (\ref{zis}) and 
(\ref{eq1}). Then 
\begin{equation}\label{ito}
g \ge G_{2,1}(\alpha).
\end{equation}
}

Proof. 
As for the$15\times15$ matrix  
$$
\hbox{\got G} =
\left(
\begin{array}{ccccccccccccccc}1 & 0 & 0 & 0 & 1 & 0 & 1 & 1 & 0 & 0 & 0 & 0 & 0 & 0 & 1
\cr
 -1 & 0 & 0 & 0 & 0 & 0 & 0 & 0 & 0 & 0 & 0 & 0 & 0 & 0 & 0
\cr 
0 & -1 & 0 & 0 & 0 & 0 & 0 & 0 & \alpha & 0 & 0 & 0 & 0 & 0 & 0
\cr
 0 & 0 & 0 & -1 & 0 & 0 & 0 & 0 & 0 & \alpha & 0 & 0 & 0 & 0 & 0
\cr
 0 & 0 & 0 & 0 & 0 & -1 & 0 & 0 & 0 & 0 & 0 & 0 & \alpha & 0 & 0
\cr 
0 & 0 & 0 & 0 & 0 & 0 & 0 & -1 & 0 & 0 & 0 & 0 & 0 & 0 & \alpha
\cr
 0 & 0 & 0 & 0 & 0 & 0 & 0 & 0 & 0 & 0 & -1 & 1 & 0 & 0 & 0
\cr
 0 & 0 & 0 & 0 & 0 & 0 & 0 & 0 & 0 & 0 & 0 & 0 & -1 & 1 & 0
\cr
 0 & 0 & -1 & 1 & 0 & 0 & 0 & 0 & 0 & 0 & 0 & 0 & 0 & 0 & 0
\cr
 0 & 0 & 0 & -1 & 1 & 0 & 0 & 0 & 0 & 0 & 0 & 0 & 0 & 0 & 0
\cr
 0 & 0 & 0 & 0 & 0 & -1 & 1 & 0 & 0 & 0 & 0 & 0 & 0 & 0 & 0
\cr 
0 & 0 & 0 & 0 & 0 & 0 & 0 & 0 & -1 & 0 & 0 & 0 & 0 & 0 & 0
\cr 
0 & 0 & 0 & 0 & 0 & 0 & 0 & 0 & 0 & g & -1 & 0 & 0 & 0 & 0
\cr 0 & 0 & 0 & 0 & 0 & 0 & 0 & 0 & 0 & 0 & 0 & g & -1 & 0 & 0
\cr
 0 & 0 & 0 & 0 & 0 & 0 & 0 & 0 & 0 & 0 & 0 & 0 & 0 & g & -1\end{array}\right).$$
  we have ${\rm det} \,\hbox{\got G} \ne 0$,
 the set
$$
{\cal C} =
\{ \hbox{\got X} \in {\cal R}':\,\,\, \text{coordinats of the point} \,\,\hbox{\got X}\,\,
\text{satisfy (\ref{zis})}\,\}
$$
is a simplicial cone in ${\cal R}'$.  We may calculate the 
 coordinates of the vertex of ${\cal C}$. If we substitute them into (\ref{eq1}), we see that $g$ will be the root of the polynomial $f_{2.1}^{[\alpha]}$.
So we see that  the condition 
$$
{\cal C} \cap \{ \hbox{\got X}: \,\,\, \xi_\nu +X_{\nu+1} =\xi_{r_1-1}+X_{r_1}\} \neq\varnothing
$$
is equivalent to (\ref{ito}).$\Box$

Now we prove (\ref{oo2}). Suppose that  $\hbox{\got i}(\Theta ) = 2$ and $\hat{\omega} \le\frac{1}{2}$. We take $\alpha <\hat{\omega}(\Theta) $ and take the 4-tiple $ (\nu,r_{1}r_{2},k)=(\nu_l,r_{l,1}r_{l,2},k_l)$ with $l$ large enough. We may suppose that (\ref{aal}) 
holds for all the indices under consideration.
Suppsose that $ q_{j+1} \le q_j ^{g}$ for  $j \in \{ \nu, r_1,r_2,k\}$.
Put
$$
\xi_j = \log \zeta_j  \,\,\,\, X_j = \log q_j  ,\,\,\,\, j = 1,2,3,... .
$$
From Proposition 1 we see that 
 all the 
inequalities (\ref{zis}) are satisfied with the only one exception. Instead of the first inequality there will be the  inequality 
$$
\xi_{\nu-1}+\xi_{r_1-1}+\xi_{k-1}+\xi_k+ x_{k+1} \ge \gamma(\Theta),
$$
 where $\gamma(\Theta)$ is bounded as $\l\to \infty$.  
At the same time instead of (\ref{eq1}) we have
$$
|\xi_\nu +X_{\nu+1}-(\xi_{r_1-1}+X_{r_1})| \le \delta (\Theta),
$$
 where $\delta(\theta)$ is bounded as $\l\to \infty$.
We take into account that  $|\log \xi_j|, |X_j| \to \infty$  for all indices under consideration, as $l \to \infty$.
So from Lemma 1 we see that for any positive $\varepsilon$ for large $l$ either
$ X_{\nu+1} \ge (G_{2,1}(\hat{\omega}) -\varepsilon) X_\nu$,
or
$ X_{r_1+1} \ge( G_{2,1}(\hat{\omega}) -\varepsilon)X_{r_1}$,
or
$ X_{r_2+1} \ge( G_{2,1}(\hat{\omega}) -\varepsilon)X_{r_2}$,
or
$ X_{k+1} \ge (G_{2,1}(\hat{\omega}) -\varepsilon) X_k$,
So there exist $j \in \{ \nu, r_1,r_2,k\}$ such that
$ q_{j+1} \ge q_j^{G_{2,2}(\hat{\omega})-\varepsilon}$ and everything follows.

{\bf 6. Case  $ \hbox{\got i}(\Theta) = 2$ and $\hat{\omega}\ge \frac{1}{2}$.}

Analogically, we consider the same vector $\hbox{\got X}$ and a similar, but distinct system of inequalities

\begin{equation}\label{zis2}
\begin{cases}
\xi_{\nu-1}+\xi_{r_2-1}+\xi_{k-1}+\xi_k+ X_{k+1} \ge 0,\cr
\xi_{r_1-1}+\xi_{r_1}+ X_{r_1+1}\ge 0,\cr
\xi_{j-1}-\alpha X_j \le 0,\,\,\, j = \nu,\nu+1,r_1+1,r_2+1, k+1\cr
X_{r_1+1}\le X_{r_2},\,\,\, X_{r_2+1} \le X_k,\cr
\xi_{r_1}\ge\xi_{r_2-1},\,\,\, \xi_{r_2}\ge \xi_{k-1},\cr
X_{j+1} \le gX_{j},\,\,\, j =\nu, r_1,r_2, k.  
\end{cases}
\end{equation}

We deal with the same equation (\ref{eq1}).

{\bf Lemma 2.}\,\,{\it Suppose that $\alpha, g >0$.
Suppose that there exists $\hbox{\got X} \in {\cal R}'$ such that its coordinates satisfy (\ref{zis2}) and 
(\ref{eq1}). Then 
\begin{equation}\label{ito}
g \ge G_{2,2}(\alpha).
\end{equation}
}
The proof is quite similar to the proof of Lemma 1.
$\Box$
\\
Now we sketch a proof for  the second part of (\ref{oo2}). The argument from the previous case is valid if we take into account that from the
 Proposition 1 it follows that
$$
\xi_{r_1-1}+\xi_{r_1}+ X_{r_1+1}\ge \gamma_1(\Theta)
$$
where  $\gamma_1(\Theta)$ is a constant. Now we may apply Lemma 1 similarly to the previous case and obtain the desired statement.

{\bf 7. Case  $ \hbox{\got i}(\Theta) = 3$: a sketch.}

Now we consider a vector
$$
\hbox{\got X} = (\xi_{\nu-1},\xi_\nu, \xi_{r_1-1},\xi_{r_1},\xi_{r_2-1},\xi_{r_2},\xi_{r_3-1},\xi_{r_3},\xi_{k-1},\xi_k, X_\nu, X_{\nu+1},
X_{r_1}, X_{r_1+1}, X_{r_2}, X_{r_2+1},  X_{r_3}, X_{r_3+1}, X_k, X_{k+1}).
$$
and a system of inequalities.
\begin{equation}\label{zis3}
\begin{cases}
\xi_{\nu-1}+\xi_{r_3-1}+\xi_{k-1}+\xi_k+ X_{k+1} \ge 0,\cr
\xi_{j-1}-\alpha X_j \le 0,\,\,\, j = \nu,\nu+1,r_2+1,r_3+1, k+1\cr
X_{r_1+1}\le X_{r_2},\,\,\, X_{r_2+1} \le  X_{r_3},  X_{r_3+1} \le X_k,\cr
\xi_{r_1}\ge\xi_{r_2-1}\ge\xi_{r_2}\xi_{r_3-1},\,\,\, \xi_{r_3}\ge \xi_{k-1},\cr
X_{j+1} \le gX_{j},\,\,\, j =\nu, r_1,r_2, r_3, k.  
\end{cases}
\end{equation}
with the equation  (\ref{eq1})

To obtain (\ref{oo3}) we need

{\bf Lemma 3.}\,\,{\it Suppose that $\alpha, g >0$.
Suppose that there exists $\hbox{\got X} \in {\cal R}'$ such that its coordinates satisfy (\ref{zis3}) and 
(\ref{eq1}). Then 
\begin{equation}\label{ito}
g \ge G_{3}(\alpha).
\end{equation}
}
To prove Lemma 3 one should consider a 19-dimensional cone in the subspace $\{ X_\nu = 0\}$ defined by (\ref{zis3}).

\end{document}